\documentclass{article}
\usepackage{graphicx} 
\usepackage{amsmath}
\usepackage{amssymb}
\usepackage{physics}
\usepackage[noblocks]{authblk}
\usepackage{hyperref}
\usepackage{wrapfig}
\title{On the Closed-Form Solution to the Interior Goat Problem in One Dimension}
\author[1]{Nathaniel Dene Hoffman\footnote{Email: \href{https://www.youtube.com/watch?v=dQw4w9WgXcQ}{dene@cmu.edu}}}
\affil[1]{Department of Recklessly Applied Mathematics,\newline University of San Serriffe}
\date{April 1, 2023}

\begin{document}

\maketitle

\begin{abstract}
    The ``interior goat problem" is a numerical puzzle which has haunted mathematicians for nearly three centuries. Closed-form solutions have only been found for the two-dimensional case, and although approximations have been established for higher-dimensional versions, the mathematical literature lacks a solution to the one-dimensional interior goat problem. To amend this oversight, we provide a closed-form solution to the interior goat problem in one dimension and speculate on continuations to even lower dimensional cases.
\end{abstract}

\section{Introduction}
Since the inception of the so-called ``exterior goat problem" in the 1749 edition of The Ladies' Diary~\cite{Heath}, there have been several efforts~\cite{Glaister1991,Fraser1982,Fraser1984,Meyerson1984,Ullisch2020} to find solutions to the related ``interior goat problem" (and various extensions) posed by Charles E. Myers in the first edition of American Mathematical Monthly~\cite{Heaton1894}, which reads as follows:

\begin{quote}
    A circle containing one acre is cut by another whose centre is on the circumference of the given circle, and the area common to both is one-half acre. Find that radius of the cutting circle.
\end{quote}

The elegant solution recently given by Ullisch~\cite{Ullisch2020} is unique in that it is a closed-form solution:
\begin{equation}
    R = 2r\cos\left(\frac{1}{2}\frac{\displaystyle\oint_{\abs{z - 3\pi/8}=\pi/4} z \dd{z} / (\sin z - z \cos z - \pi/2)}{\displaystyle\oint_{\abs{z - 3\pi/8}=\pi/4} \dd{z} / (\sin z - z \cos z - \pi/2)}\right)
\end{equation}
where $r$ is the radius of the circular fence and $R$ is the length of the tether. In the case where $r=1$, Ullisch approximates the solution as $R \approx 1.15872847302$.

\begin{figure}[h]
\centering
\includegraphics[width=0.5\textwidth]{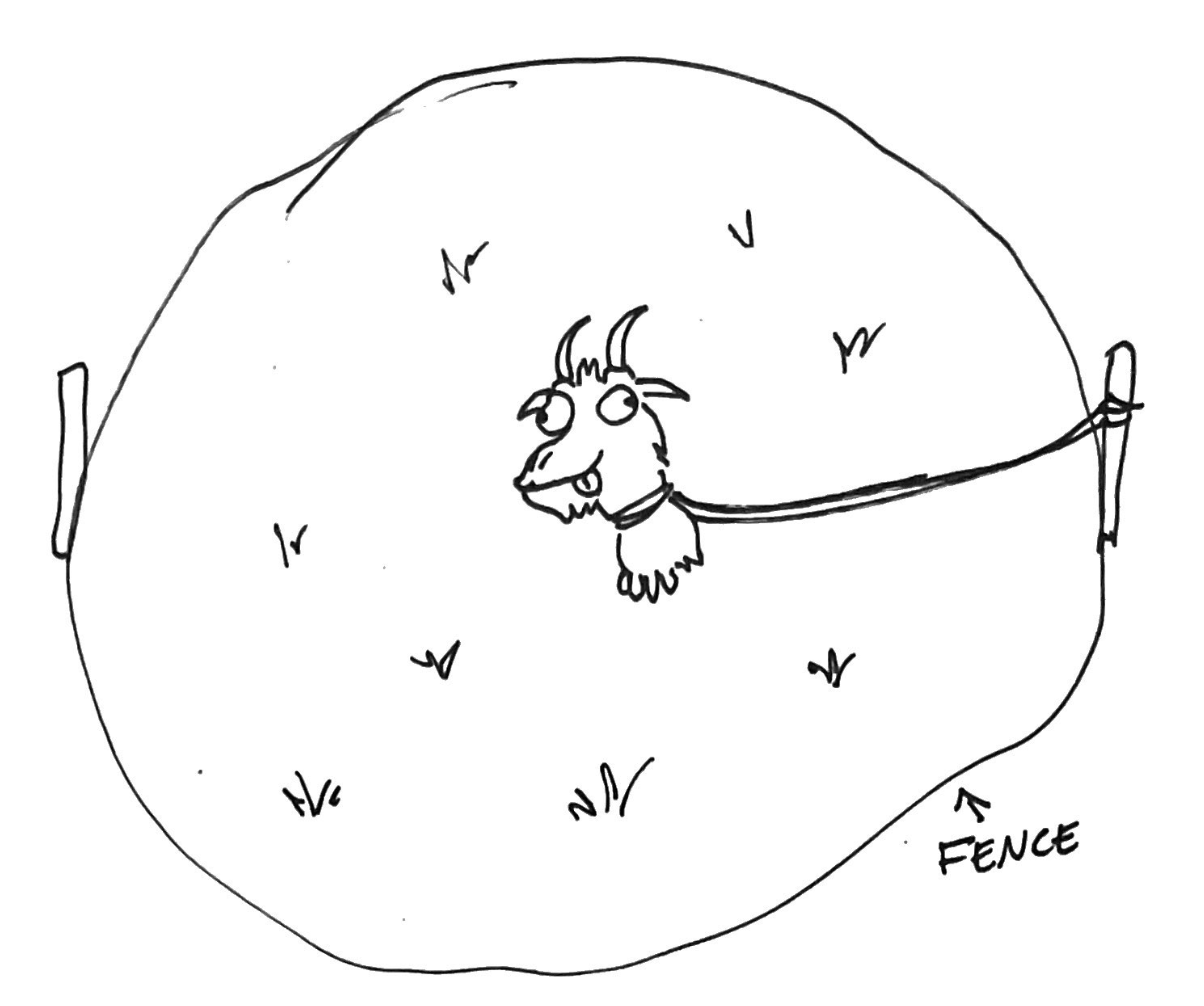}
\caption{The interior goat problem in two dimensions. This goat is cruelly tethered so that he may only graze half the field.}
\label{fig:2D}
\end{figure}

In 1984, Fraser~\cite{Fraser1984} found general solutions for extensions to $n$-dimensional fields (with $n+1$-dimensional goats, of course), which happened to be correct aside from a slight error in a proof, corrected by Meyerson~\cite{Meyerson1984}. Fraser gives $R_n = k_n r$ where $k_n = 2\cos\beta$. Furthermore, $\beta$ is defined as the angle between the line connecting the center of the field with the tether's non-goat endpoint to the furthest point along the fence which the goat can reach. Approximate values of $\beta$ (and furthermore $k_n$) can be obtained by solving the following equation numerically:
\begin{equation}
    (2\cos\beta)^n \int_{(\pi/2) - \beta}^{\pi/2}\cos^n\theta \dd{\theta} = \int_0^{2\beta - (\pi/2)}\cos^n\theta \dd{\theta}
\end{equation}

Fraser gives approximate solutions for the cases $n=2,3,4,5$, where the $n=2$ case is simply
\begin{equation}
    \sin 2\beta - 2\beta\cos 2\beta = \frac{\pi}{2} \implies k_2 \approx 1.158728,
\end{equation}
the same numerical answer found in the closed-form solution by Ullisch. Additionally, \cite{Fraser1984} and \cite{Meyerson1984} both come to the conclusion that $\lim_{n\to\infty} k_n = \sqrt{2}$.

However, it appears that the instance where $n=1$ is strikingly absent from the related literature. In the following section, we propose a closed-form solution for this important case.

\begin{figure}[h]
\centering
\includegraphics[width=0.7\textwidth]{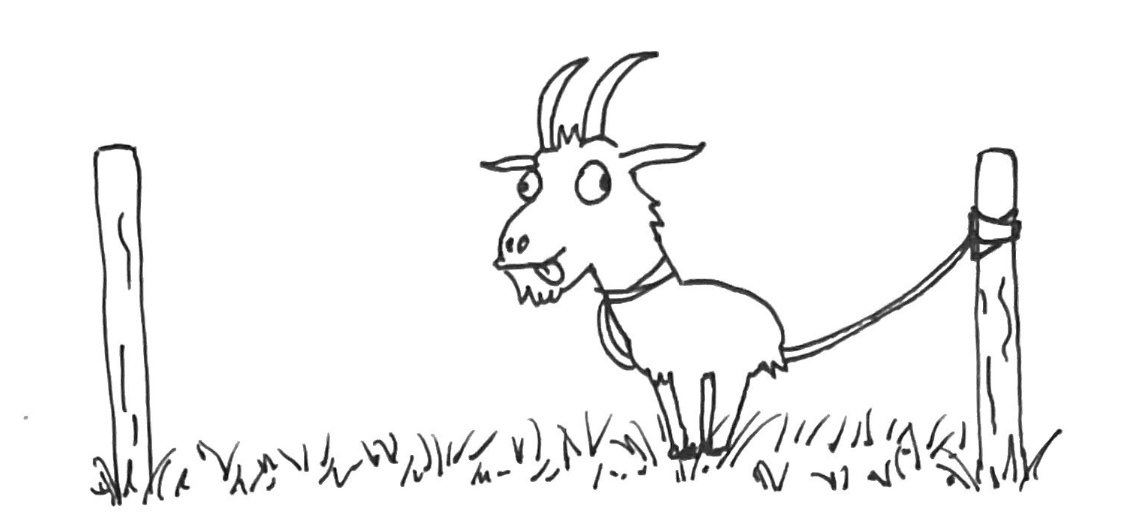}
\caption{The interior goat problem in one dimension. The tether length required to restrict him to half the field was previously undefined in the mathematical literature.}
\label{fig:1D}
\end{figure}

\section{The Closed-Form Solution}
We begin with the work of Fraser, now setting $n=1$:
\begin{align}
    (2\cos\beta)^1 \int_{(\pi/2) - \beta}^{\pi/2}\cos^1\theta \dd{\theta} &= \int_0^{2\beta - (\pi/2)}\cos^1\theta \dd{\theta} \notag\\
    2\cos\beta\left[\eval{\sin\theta}_{\theta=(\pi/2) - \beta}^{\pi/2}\right] &= \left[\eval{\sin\theta}_{\theta=0}^{2\beta - (\pi/2)}\right] \notag\\
    2\cos\beta\left[\sin(\frac{\pi}{2}) - \sin(\frac{\pi}{2} - \beta)\right] &= \left[\sin(2\beta - \frac{\pi}{2}) - \sin(0)\right] \notag\\
    2\cos\beta\left[1 - \sin(\frac{\pi}{2} - \beta)\right] &= \left[\sin(2\beta - \frac{\pi}{2}) - 0\right]
\end{align}
This might seem quite daunting for the uninitiated, but we can use a clever trigonometric identity to further reduce the equation:
\begin{equation}\label{eq:clever_trig_identity}
    \sin(\alpha - \beta) = \sin(\alpha)\cos(\beta) - \cos(\alpha)\sin(\beta)
\end{equation}
Therefore, by applying Equation~\ref{eq:clever_trig_identity}, we find
\begin{align}
    2\cos\beta\left[1 - 1\cdot\cos(\beta) - 0\cdot\sin(\beta)\right] &= \left[\sin(2\beta)\cdot 0 - \cos(2\beta)\cdot 1\right] \notag\\
    2\cos\beta\left[1 - \cos(\beta)\right] &= -\cos(2\beta) \notag\\
    2\cos(\beta) - 2\cos[2](\beta) + \cos(2\beta) &= 0 \notag\\
    2\cos(\beta) - 1 - \cos(2\beta) + \cos(2\beta) &= 0 \notag\\
    \cos(\beta) &= \frac{1}{2} \notag\\
    \implies \beta &= \pm\frac{\pi}{3} + 2\pi m,\quad m\in\mathbb{Z}
\end{align}
Finally, we can use this $\beta$ to find the solution:
\begin{equation}
    k_1 = 2\cos\beta = 2\cos(\frac{\pi}{3}) = 2\cdot \frac{1}{2} = 1
\end{equation}
where we ignore the $\pm$ and the $+2\pi m$ since the cosine function is even and has a period of $2\pi$.

\section{Discussion}
How should we interpret this solution? First, what do we mean by a one-dimensional pasture? The higher-dimensional pastures are defined by an $n$-ball enclosed by an $n$-sphere fence. Generally, the volume of an $n$-ball is given by
\begin{equation}
    V_n(r) = \frac{\pi^{n/2}r^n}{\Gamma\left(\frac{n}{2} + 1\right)}
\end{equation}
so
\begin{equation}
    V_1(r) = \frac{\pi^{1/2}r^1}{\Gamma\left(\frac{3}{2}\right)} = 2r
\end{equation}
To graze half of the volume of the pasture, we must therefore allow the goat to have access to exactly $\frac{1}{2}\cdot 2r = r$ units of the pasture. From the solution for $k_1$ given above, we can then conclude that the tether length $R = k_n r$ must be $R = r$, the radius of the field. This certainly feels right\textemdash I've asked a bunch of people, and they usually agree with me.

\begin{wrapfigure}{r}{0.35\textwidth}
\centering
\includegraphics[width=0.9\linewidth]{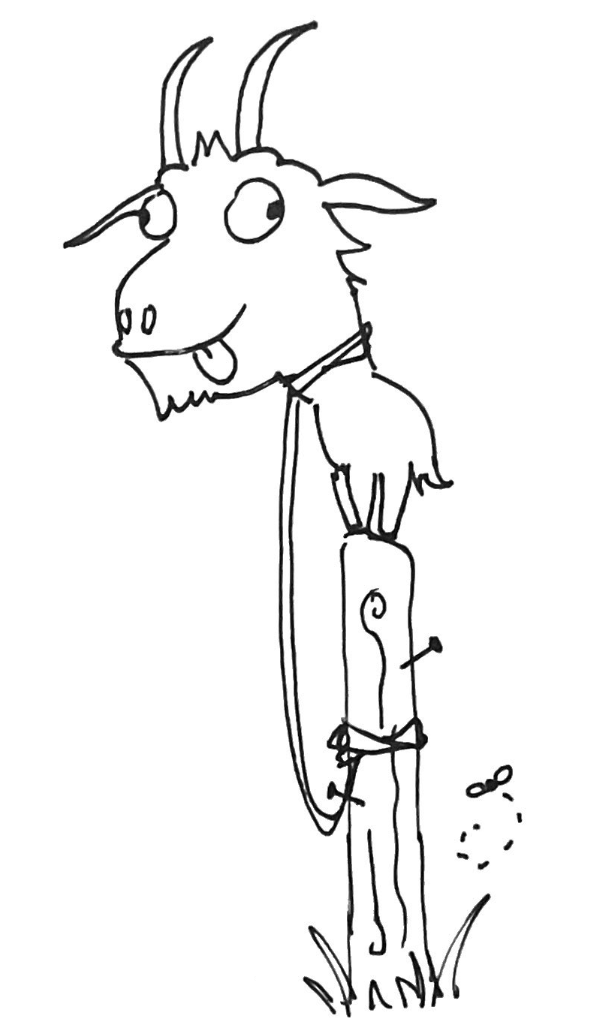}
\caption{The interior goat problem in zero dimensions. This goat can graze half the field as well as the whole field!}
\label{fig:0D}
\end{wrapfigure}

How does one fence in a one-dimensional pasture? Since we define the fence as an $n$-sphere, the pasture will only have two fence posts, one at either end of the field, which will surely save a fortune on construction costs. Furthermore, for the $n=0$ case, we have
\begin{equation}
    \beta = 2\beta - \frac{\pi}{2} \implies \beta = -\frac{\pi}{2} \implies k_0 = 0
\end{equation}
whereas, of course,
\begin{equation}
    V_0(r) = \frac{\pi^0 r^0}{\Gamma(1)} = 1
\end{equation}
so a one-dimensional goat in a zero-dimensional field can graze half the area while standing still (see Figure~\ref{fig:0D}).

\begin{figure}[h]
\centering
\includegraphics[width=0.4\textwidth]{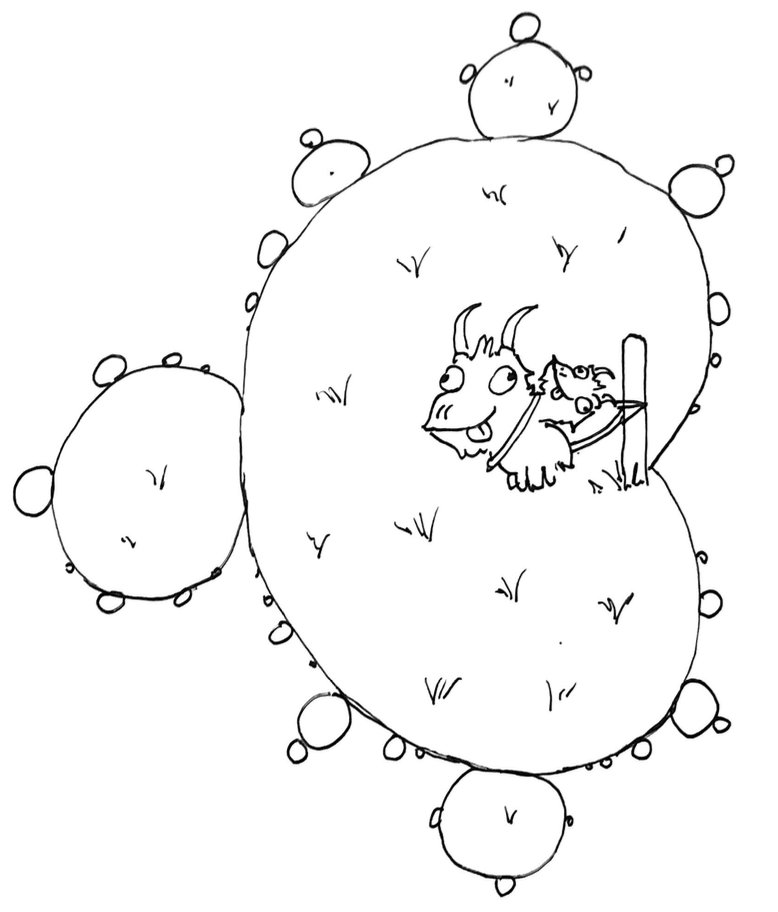}
\caption{A goat in a fractal field. This problem seems pretty difficult, so I'll leave it as an exercise for the reader.}
\label{fig:fractal}
\end{figure}

We present this result to clear up any misconceptions surrounding the solutions to the lower-dimension goat problems. As we can clearly see, the solutions defy na\"ive expectations. Future studies in fractional goat dimensions may prove useful (see Figure~\ref{fig:fractal}), and the author would also propose an \textit{in vivo} study for the $n=1$ case, given such a two-dimensional goat can be found\footnote{It is currently unknown whether there exists a two-dimensional species of goat whose body isn't completely split in half by its digestive tract.} and funding can be acquired to build a one-dimensional pasture.

\section{Acknowledgements}
The author would like to thank Geoffrey Chaucer, Eloy d'Amerval, and Eduard de Dene (no relation) for their helpful contributions. The author would also like to thank AJ Hoffman (yes relation) for his illustrations. Additionally, the author thanks Ingo Ullisch for his creative solution to the original problem which inspired this paper.

\bibliography{goat}{}
\bibliographystyle{plain}

%
%
%
%
%
%

\end{document}